\documentclass[12pt]{amsart}
\usepackage{amsfonts, amsmath, amssymb, amscd, latexsym,graphicx}

\newtheorem{thm}{Theorem}[section]
\newtheorem{prop}[thm]{Proposition}
\newtheorem{lem}[thm]{Lemma}
\newtheorem{cor}[thm]{Corollary}

\begin{document}

\title{Bi-orderings on pure braided Thompson's groups}

\author{Jos\'e Burillo}
\address{
Departament de Matem\'atica Aplicada IV,  Escola Polit\`{e}cnica
Superior de Castelldefels, Universitat Polit\'ecnica de Catalunya,
Av. del Canal Ol\'{\i}mpic s/n, 08860 Castelldefels, Barcelona, Spain.}
\email{burillo@ma4.upc.edu}
\thanks{The first author acknowledges support from
grant MTM2005-04104}

\author{Juan Gonz\'alez--Meneses}
\address{Departamento de \'Algebra, Facultad de Matem\'aticas, Universidad de Sevilla, Apdo. 1160, 41080 Sevilla, Spain.}
\email{meneses@us.es}
\thanks{The second author was partially supported by MTM2004-07203-C02-01 and FEDER}

\date{August 23, 2006}


\begin{abstract}
In this paper it is proved that the pure braided Thompson's group $BF$ admits a bi-order, analog to the bi-order of the pure braid groups.
\end{abstract}

\maketitle

\section*{Introduction}

Braid groups have been a constant object of study all along the 20th
century, and have exerted a strong fascination on algebraists since
their inception, due to their rich and interesting properties and
useful applications to other branches of mathematics. Literature
about braid groups is abundant; we will mention only a small sample
of it, namely, the seminal works of Artin~\cite{Artin1925,Artin}, a
more modern text~\cite{Hansen}, and a recent
survey~\cite{Birman-Brendle}. In particular, the properties which
will be of interest here are their orderings and bi-orderings,
see~\cite{dehornoy}, \cite{5authors} and \cite{KR} for details.
Orderings of groups are also a classical subject of study, related,
for instance, to the existence of zero divisors in the group ring.

Thompson's groups have been studied since the late 1960s, initially
as examples of infinite, finitely presented simple groups ---$V$
being the first such a group known---, but in the subsequent years
they showed other, equally striking properties, for instance, $F$
was the first torsion-free $FP_{\infty}$ group, see
\cite{BG}. For details and many proofs of the properties of
Thompson's groups see the excellent introduction in \cite{CFP}.

Recently and independently, Brin \cite{brin} and Dehornoy
\cite{dehornoyBV} have introduced braided versions of Thompson's
groups, which show a mixture of properties of both. They are a very
natural extension of Thompson's group $V$, where the permutations,
characteristic in $V$, have been replaced by braids. Hence, the
braided Thompson's group $BV$ (in Brin's notation) is a torsion-free
version of $V$, sharing many of its properties, for instance, finite
presentation, with a presentation very similar to that of $V$. One
can think of the braided Thompson's groups as ``Artin" versions of
Thompson's groups $V$ and $F$, very much in the same way as braid
groups are the corresponding Artin groups for the permutation
groups.

The main concern of this paper are some properties that braided
Thompson's groups inherit from braid groups, namely, orderings and
bi-orderings. It is a well-known fact that braid groups are
orderable but not bi-orderable, and pure braid groups are
bi-orderable. The same thing happens in the braided Thompson's
groups: the braided version of $V$, analog to $B_\infty$, is
orderable but not bi-orderable, and its pure subgroup $BF$ is
bi-orderable, as it models on $P_\infty$. This last statement is the
main result of this paper, Theorem~\ref{T:BF bi-orderable}. On the way
to proving this result, we study a somewhat non-standard version of
the pure braid group on infinitely many strands, where there are $n$
different embeddings of $P_n$ into $P_{n+1}$, each one defined by
splitting a given strand in two. The direct limit obtained this way
is an infinitely generated group, already appeared in \cite{BBCS},
and which is also bi-orderable, a fact whose proof is the crucial
ingredient of the proof of the main result, and takes up the largest
and most technical part of Section 3. The previous sections of the
paper are dedicated to set up the table for the main result, the
first section with an introduction to the braided Thompson's groups,
and the second one about orderings and the ordered groups that will
be involved in the construction of the bi-order of $BF$.

The first author would like to acknowledge the great hospitality of
the Departamento de \'Algebra of the Universidad de Sevilla, during
a visit to which this work was developed.

\section{The braided Thompson's groups}

Thompson's group $F$ is the group of piecewise-linear,
orientation-preserving homeomorphisms of the interval $[0,1]$
whose breakpoints are dyadic integers, and whose slopes are powers
of 2. An element of $F$ can be understood as a map between two
subdivisions of the unit interval, subdivisions into intervals of
lengths $1/2^n$ with dyadic endpoints, where the subintervals are
mapped linearly in order-preserving fashion. See \cite{CFP} for an
excellent introduction to Thompson's groups and proofs of their
basic properties.

\begin{figure}[ht]
\centerline{\includegraphics{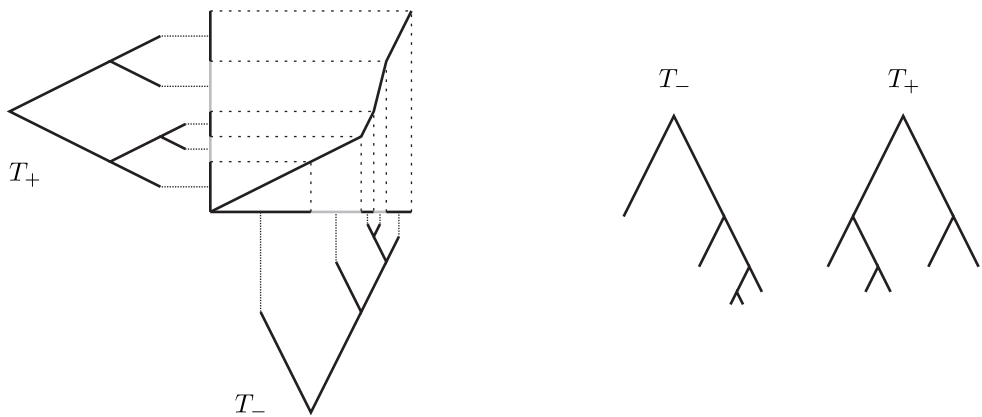}} \caption{An element of
$F$, and its reduced pair of binary trees $(T_-,T_+)$.}
\label{felem2}
\end{figure}

Such a subdivision of the unit interval is in correspondence with a
binary tree, in the usual way. Hence, an element of $F$ can be
represented as a pair of binary trees $(T_-,T_+)$, with the same
number of leaves, as is shown in Figure~\ref{felem2}. An element
admits more than one representation as a pair of binary trees by
further subdividing an interval and its target, but each element
admits a unique reduced element whose subdivision is optimal. See
again \cite{CFP} for details.

\begin{figure}[ht]
\centerline{\includegraphics{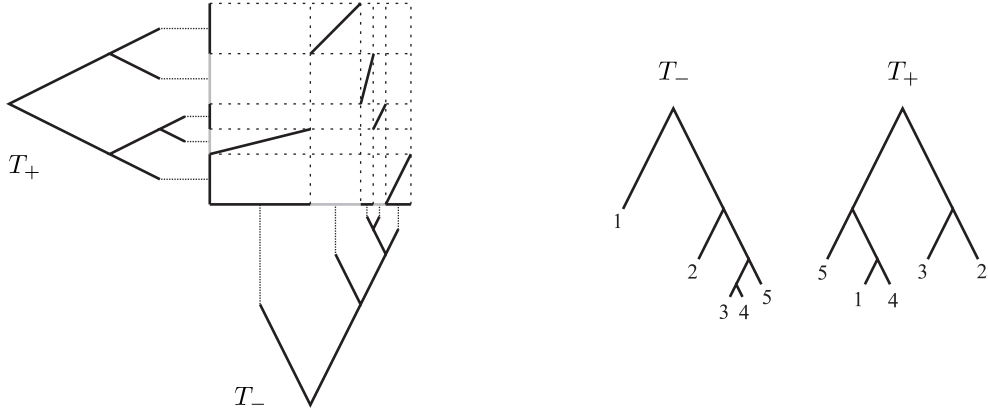}} \caption{An element of
$V$, and its pair of binary trees $(T_-,T_+)$. The labels of the
leaves describe the permutation $\pi$.} \label{velem2}
\end{figure}

Thompson's group $V$ is represented in a similar way by subdivisions
of the interval, but now the subintervals are still mapped linearly to each
other, but they can be permuted, thus not necessarily preserving their order. So
an element of $V$ is represented as a triple $(T_-,\pi,T_+)$, where
the binary trees $T_-$ and $T_+$ have $n$ leaves, and $\pi$ is a
permutation in $S_n$ which indicates how the leaves are mapped to
each other. See \cite{CFP} \and Figure~\ref{velem2} to understand
this interpretation of elements of $V$.

The group $V$ is infinite, simple and finitely presented, being the
first such group historically known. Observe that fixing a tree $T$,
and varying $\pi$, the elements $(T,\pi,T)$ form a subgroup of $V$
isomorphic to $S_n$. Hence, $V$ has torsion, and indeed, it contains
all finite groups as subgroups. As a contrast, $F$ is torsion-free.

\begin{figure}[ht]
\centerline{\includegraphics{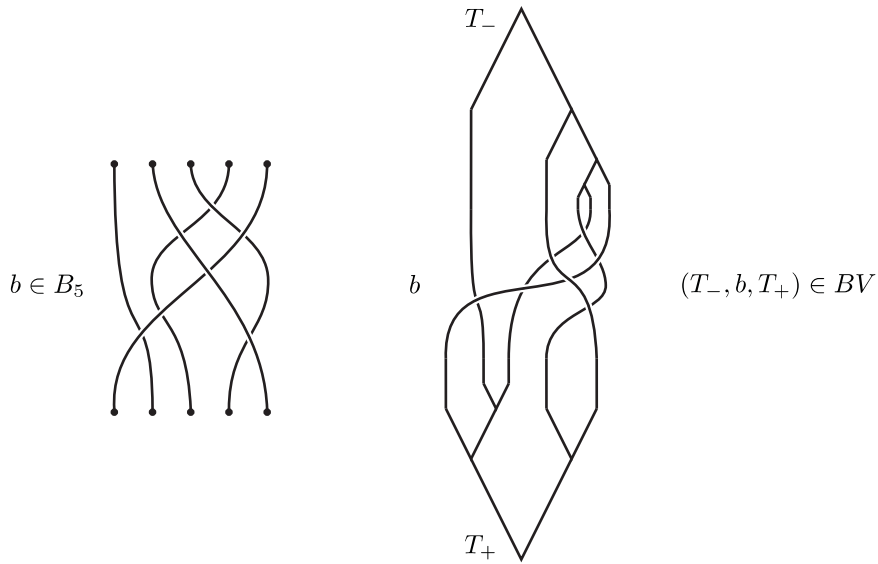}} \caption{An element
of $B_5$ and an element of $BV$.} \label{belem_bvelem}
\end{figure}

We will denote Artin's braid group in $n$ strands by $B_n$. It is
the fundamental group of the configuration space of $n$ unordered
points in the plane. See~\cite{Birman} for details. Its elements are
usually visualized as a set of $n$ disjoint strands in 3-space,
as in the left hand side of Figure~\ref{belem_bvelem}, and it admits
the following well-known presentation~\cite{Artin}:
\begin{equation}
\label{E:braid group presentation} B_n= \left< \sigma_1,\ldots,
\sigma_{n-1} \left|
\begin{array}{ll}
 \sigma_i \sigma_j = \sigma_j \sigma_i  & \mbox{if } \; |i-j|>1,  \\
 \sigma_i \sigma_j \sigma_i = \sigma_j \sigma_i \sigma_j &  \mbox{if }\; |i-j|=1.\ \end{array}\right.
 \right>.
\end{equation}
Let
$$
\rho_n:B_n\longrightarrow S_n
$$
be the homomorphism from the braid group onto the symmetric group,
which maps a braid to the permutation it induces on the base points.
The kernel of $\rho_n$ is $P_n$, the subgroup of pure braids in $n$
strands.

The braided group $BV$ (which means ``braided $V$") is a
torsion-free version of $V$ using braids instead of permutations.
Elements of $BV$ are seen as triples $(T_-,b,T_+)$, where the trees
$T_-$ and $T_+$ have $n$ leaves, and $b\in B_n$ is a braid in $n$
strands. The braid is understood as joining the leaves of $T_-$ with
those of $T_+$, see the right hand side of Figure~\ref{belem_bvelem}
for an example. An element of $BV$, as it happens in other
Thompson's groups, admits many representatives, by adding carets to
the trees and splitting the corresponding strands. Hence, to
multiply two elements in $BV$ we only need to subdivide the trees
until we find representatives with matching trees: given two
elements $(T_-,b,T_+)$ and $(T'_-,b',T'_+)$, we construct, by adding
carets to the trees and splitting strands into parallel ones, two
representatives $(\overline T_-,\overline b,\overline T_+)$ and
$(\overline T'_-,\overline b',\overline T'_+)$ such that $\overline
T_+=\overline T'_-$, and then the product is the triple $(\overline
T_-,\overline b\overline b',\overline T'_+)$, where the product of
the two braids takes place in the corresponding braid group, with as
many strands as leaves in the trees. Observe that since $\overline
T_+=\overline T'_-$, the two braids have the same number of strands.
For details on $BV$, see \cite{brin}, \cite{dehornoyBV} and
\cite{BBCS}.

As it happens in $B_n$, forgetting the braid and focusing on the
corresponding permutation gives a homomorphism
$$
\widetilde\rho :BV\longrightarrow V
$$
such that
$$
\widetilde \rho (T_-,b,T_+)=(T_-,\rho(b),T_+).
$$
Notice that the element in Figure~\ref{belem_bvelem} maps to the
element in Figure~\ref{velem2}.

\begin{figure}[ht]
\centerline{\includegraphics{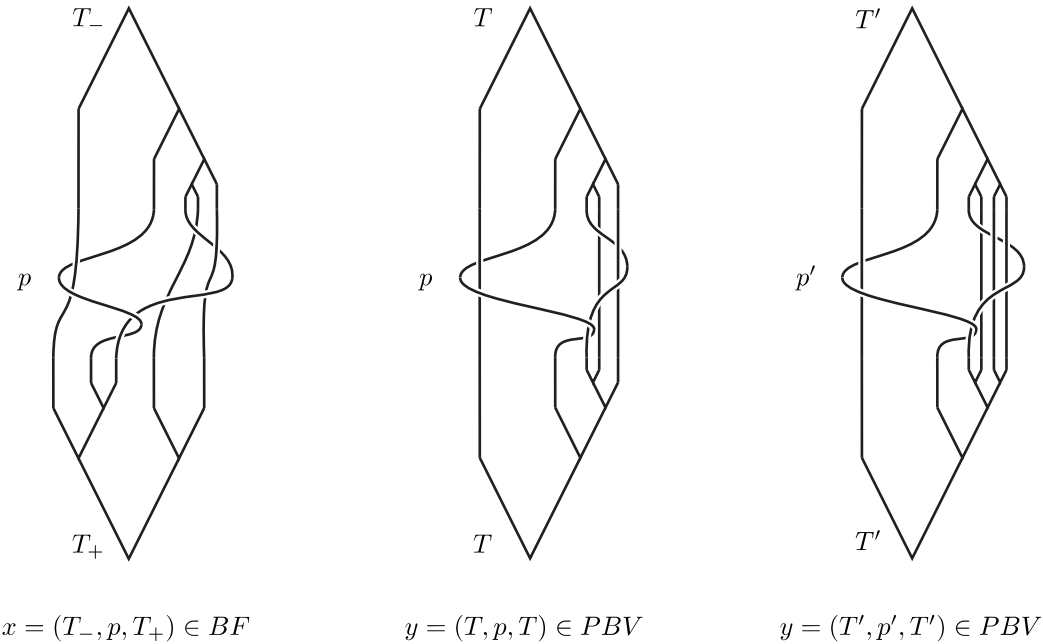}} \caption{An
element $x\in BF$, and an element $y\in PBV$ represented in two
different ways.} \label{bfelem_pbvelem}
\end{figure}

 Finally, the group $BF$ is the subgroup of $BV$ of those
elements $(T_-,p,T_+)$ whose braid is pure. See an example in the
left hand side of Figure~\ref{bfelem_pbvelem}. Observe that if $p$
is pure, then $\rho(p)$ is the identity permutation, so then
$(T_-,\rho(p),T_+)$ is actually an element of $F$ inside $V$. Then,
it is clear that $BF={\widetilde \rho^{-1}}(F)$, hence the name $BF$
(``braided $F$"). For details on $BF$, in particular for a finite
presentation, see \cite{BBCS}.

A group which will be of great relevance later is the group
$PBV=\ker\widetilde\rho$. Clearly, since  $(T_-,p,T_+)$ is in $PBV$
if $\widetilde\rho(T_-,p,T_+)=1$, we must have that $p$ is a pure
braid, and also $T_-=T_+$. Then, $PBV$ is the subgroup of $BV$
(actually inside $BF$) of all those elements where the two trees are
the same and the braid is pure. See an example in
Figure~\ref{bfelem_pbvelem}. Observe that if an element has a
representative where the two trees are the same and the braid is
pure, then all representatives satisfy these two conditions.

\begin{figure}[ht]
\centerline{\includegraphics{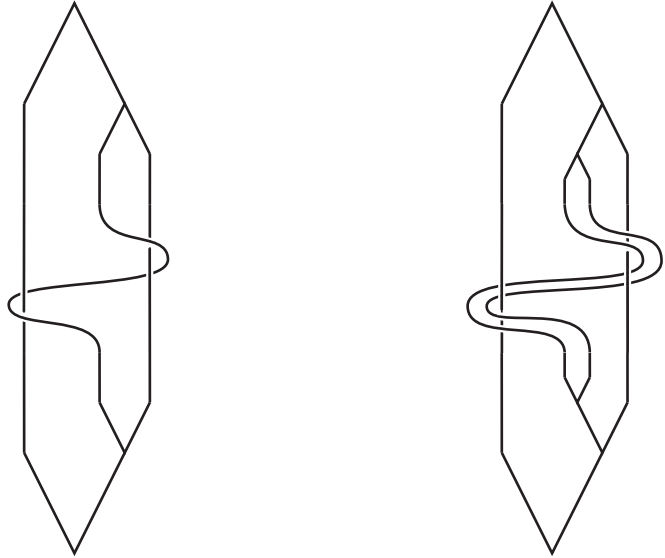}} \caption{An example of
an element in $P_{3,T}$ and its image under $\alpha_{3,T,2}$, by
splitting the second strand.} \label{alphanTi}
\end{figure}

Given a tree $T$ with $n$ leaves, the subgroup of $PBV$ of the
elements represented by $(T,p,T)$, varying $p$, is a group
isomorphic to $P_n$. We will denote this particular copy of $P_n$
inside $PBV$ by $P_{n,T}$, and there are as many such subgroups
isomorphic to $P_n$ inside $PBV$ as there are trees with $n$ leaves,
i.e., the Catalan number
$$
\text{Cat}(n-1)=\frac 1n{{2n-2}\choose{n-1}}.
$$
If we now consider the tree $T'$ obtained from $T$ by attaching a
caret to the $i$-th leaf of $T$, we get another representative for
$(T,p,T)$, namely, $(T',p',T')$, where $p'$ has been obtained from
$p$ by splitting the $i$-th strand in two parallel ones (see
examples in Figure~\ref{bfelem_pbvelem} and Figure~\ref{alphanTi}).
We have then a one-to-one homomorphism
$$
\alpha_{n,T,i}:P_{n,T}\longrightarrow P_{n+1,T'}
$$
obtained via this process. The group $P_{n,T}$ can be identified
with the subgroup of $P_{n+1,T'}$ of those elements whose $i$-th and
$(i+1)$-st strands are parallel. Observe that both an element of
$P_{n,T}$ and its image under $\alpha_{n,T,i}$ represent the same
element in $PBV$. With these subgroups and maps, the following
proposition is straightforward.

\begin{prop}
The groups $P_{n,T}$, for $n>1$ and $T$ a tree with $n$ leaves,
together with the maps $\alpha_{n,T,i}$, for $1\le i\le n$, form a
direct system of groups and homomorphisms, whose direct limit is
isomorphic to $PBV$. That is,
$$
  PBV=  \lim_{\substack{\displaystyle \longrightarrow \\ T}}
  P_{n,T}.
$$
\end{prop}

In the next section we will recall the notion of orderable groups.
We will also recall that the pure braid group $P_n$ is bi-orderable,
and we will review the bi-order for the Thompson group $F$.
Showing that there is a bi-order of $P_n$ which is consistent with
the above direct system, we will be able to show that $PBV$ is
bi-orderable. This fact, together with the bi-order of $F$, will
allow us to show that $BF$ is also bi-orderable.

\section{Orderings}

\subsection{Concepts and basic properties}

A group $G$ is said to be {\em left-orderable} if there exists a
total order on its elements which is invariant under left
multiplication, that is, $a<b$ implies $ca<cb$ for every $a,b,c\in
G$ (see~\cite{MR}). Such an order is determined by the set of
positive elements, $P=\{x\in G;\ 1<x\}$, since $a<b$ if and only if
$1<a^{-1}b$, in other words, if $a^{-1}b\in P$.

Notice that every subset $P\in G$ determines a binary relation on
the elements of $G$ in the above way ($a<b \Leftrightarrow
a^{-1}b\in P$). This relation is clearly invariant under
left-multiplication. Moreover, the relation is transitive if and
only if $P$ is a semigroup, and it is antisymmetric and total if and
only if $G=P\sqcup \{1\} \sqcup P^{-1}$. Hence, a group $G$ is
left-orderable if and only if it contains a subsemigroup $P\subset
G$ such that $G=P\sqcup \{1\} \sqcup P^{-1}$.

A very simple example of left-orderable group is $\mathbb Z^n$, with
the lexicographical order. A less obvious example is
$B_n$~\cite{dehornoy,5authors}.  Notice that, in a left-orderable
group, all powers of a positive element are positive, hence
left-orderable groups are torsion-free. This shows in particular
that finite groups, as well as Thompson's group $V$, are not
left-orderable.

If a group is left-orderable, one can use the semigroup $P$ to
define a {\em right-order} $\prec$, that is a total order of its
elements which is invariant under right multiplication: we just say
that $a\prec b$ if $b\: a^{-1}\in P$. Hence, a group is
left-orderable if and only if it is right-orderable, but the two
orderings do not necessarily coincide.

A group is said to be {\em bi-orderable} if it admits a left-order
which is also a right-order. Notice that a group $G$ is bi-orderable
if and only if it admits a subsemigroup $P$ such that $G=P\sqcup
\{1\}\sqcup P^{-1}$ (thus $G$ is left-orderable), and furthermore
$P$ is closed under conjugation in G (thus the left-order is also a
right-order). The group $\mathbb Z^n$ is hence bi-orderable, with
the lexicographical order. The free group $F_n$ of rank $n$ is also
bi-orderable, as we shall see later.

Notice that, in a bi-orderable group, every conjugate of a positive
element is positive. Since products of positive elements are
positive, this implies that a bi-orderable group cannot have {\it
generalized torsion}, which means that the product of nontrivial
conjugate elements can never be trivial.  The braid group $B_n$
($n>2$) is an example of a left-orderable group which is not
bi-orderable, since it has generalized torsion. Namely, in $B_n$ one
has $\sigma_1\sigma_2^{-1} \cdot (\sigma_1 \sigma_2\sigma_1)^{-1}
\sigma_1 \sigma_2^{-1} (\sigma_1 \sigma_2 \sigma_1) = 1$. But $B_n$
contains a finite index subgroup which is bi-orderable, namely the
pure braid group $P_n$. Later in this section we will recall the
explicit bi-order that can be defined in $P_n$.

There is a well-known result (see for instance~\cite{KR}) that shows
how left and bi-orderability behave under extensions.

\begin{prop}\label{P:orders and extensions}
Consider the following exact sequence of groups:
$$
   1 \longrightarrow A \stackrel{\alpha}{\longrightarrow} B \stackrel{\beta}{\longrightarrow} C
   \longrightarrow 1.
$$
Suppose that $A$ and $C$ are left-orderable, and let $P_A$ and $P_C$
be their corresponding semigroups of positive elements. Then $B$ is
also left-orderable, an explicit left-order being defined by the
semigroup $\alpha(P_A)\cup \beta^{-1}(P_C)$.

Moreover, suppose that $P_A$ and $P_C$ define bi-orders on $A$ and
$C$ respectively. Suppose also that the above sequence splits, so
$B=C\ltimes A$. If the action of $C$ on $A$ preserves the order in
$A$ (i.e. preserves $P_A$), then $B$ is bi-orderable, an explicit
bi-order being defined by the semigroup $\alpha(P_A)\cup
\beta^{-1}(P_C)$.
\end{prop}

The above left-order (say bi-order) of $B$ can be explained as
follows: An element $b\in B$ is positive if and only if $\beta(b)$
is positive in $C$, or $\beta(b)=1$ and $b$ is positive in $A$.  In
the case of the semi-direct product, we can just say that the order
in $B=C\ltimes A$ is the lexicographical order.

\subsection{Explicit bi-orders in some groups} We will now explain
some specific bi-orders in three groups which are important for our
purposes. They are Thompson's group $F$, the free group of rank $n$,
$F_n$, and the pure braid group on $n$ strands $P_n$. The three
groups are already known to be bi-orderable, and we shall explain
the known bi-orders of $F_n$ and $P_n$, together with a review of
the ordering of $F$.

The bi-order in $F$ is defined as follows. Recall that
an element $f\in F$ is a piecewise-linear, orientation preserving
homeomorphism of the interval $[0,1]$ and all slopes are powers of 2.
An element is then positive if its first slope different from 1 is a
positive power of 2. It is not difficult to prove that it is a bi-order,
see \cite{dehornoyBV} for details and other descriptions of this bi-order.

We will now recall the usual bi-order of the free group $F_n$. A
detailed proof can be found in~\cite{KR}. It is based on the so
called {\em Magnus expansion} of the free group~\cite{MKS}, so it is
usually called the {\em Magnus ordering}. Let $F_n$ be the free
group of rank $n$, freely generated by $x_1,\ldots,x_n$. Let
$\mathbb Z[[\mathbf X]]=\mathbb Z[[X_1,\ldots,X_n]]$ be the ring of
formal series on $n$ non-commutative variables $X_1,\ldots,X_n$,
with coefficients in $\mathbb Z$. The Magnus expansion is the
homomorphism $ \varphi: \: F_n \longrightarrow \mathbb Z [[\mathbf
X]] $ defined by $\varphi(x_i)=1+X_i$, for $i=1,\ldots,n$.  Notice
that one then has $\varphi(x_i^{-1})=1-X_i+X_i^2-X_i^3+\cdots$.  It
is shown in~\cite{MKS} that $\varphi$ is injective. The elements in
Im$(\varphi)$ are formal series $f(\mathbf X)$ such that $f(\mathbf
0)=1$.

Notice that the set of monomials of $\mathbb Z [[\mathbf X]]$ can be
totally ordered: First, we order the variables by $X_1<X_2<\cdots <
X_n$. Then, given two monomials, the smallest one will be the one of
smallest degree, or in case their degrees coincide, the smallest one
in lexicographical order. One can then say that
a nontrivial element $x\in F_n$ is positive if the coefficient of
the smallest nontrivial term of $\varphi(x)-1$ is positive. This
defines a set $P$ of positive elements in $F_n$ which yields a
bi-order of $F_n$.

There is an important property of the Magnus expansion that we will
need later. Define the map $\delta: \: F_n \rightarrow \mathbb
Z[[\mathbf X]]$ to be the map (but not homomorphism) that sends the
trivial element to $0$, and every nontrivial $x\in F_n$ to the
nontrivial homogeneous form of smallest degree in $\varphi(x)-1$.
The form $\delta(x)$ is called the {\it deviation} of $x$. Notice
that $x$ is positive if and only if the coefficient of the smallest
term in $\delta(x)$ is positive. Consider now the lower central
series of $F_n$, $F_n=G_1\supset G_2 \supset G_3 \supset \cdots $.
It is shown in~\cite{MKS} that $\delta(x)$ has degree $d$ if and
only if $x\in G_d$ and $x\notin G_{d+1}$. Moreover, if we denote
$\mathbb Z[[\mathbf X]]_d$ the set of homogeneous forms of degree
$d$ in $\mathbb Z[[\mathbf X]]$, one has:

\begin{thm}\cite{MKS}\label{T:deviation}
For every $d>0$, the map $\delta$ determines a one-to-one
homomorphism from the abelian group $G_d/G_{d+1}$, under group
multiplication, to $\mathbb Z[[\mathbf X]]_d$, under addition.
\end{thm}

We end this section with a exposition of a bi-order of the pure
braid group $P_n$, given in~\cite{KR}. It comes from the above
bi-orders of free groups, together with Proposition~\ref{P:orders
and extensions}, since there is a well-known split exact sequence of
groups:
$$
 1 \longrightarrow F_{n-1} \longrightarrow P_n
 \stackrel{\eta}{\longrightarrow} P_{n-1} \longrightarrow 1,
$$
where for every braid $p\in P_n$, $\eta(p)$ is the braid obtained
from $p$ by deleting its first strand. Indeed, if all strands except
the first one are trivial in a pure braid $p$, these trivial strands
can be considered as being punctures of the plane, and the first
strand of $p$ can be considered as describing a loop in the
$(n-1)$-punctured plane. Hence $\ker\eta$ is isomorphic to the
fundamental group of the $(n-1)$-punctured plane, which is a free
group on $n-1$ generators. The above sequence clearly splits, by adding a strand to the left of a braid in $P_{n-1}$, so
$P_n=P_{n-1}\ltimes F_{n-1}$, and the action of $P_{n-1}$ on
$F_{n-1}$ preserves the Magnus ordering of $F_{n-1}$~\cite{KR}. It
follows by recurrence on $n$ that $P_n=(\cdots((F_1\ltimes F_2)
\ltimes F_3) \ltimes \cdots F_{n-2}) \ltimes F_{n-1}$, and $P_n$ is
bi-orderable by Proposition~\ref{P:orders and extensions}, an
explicit bi-order being the lexicographical order in
$(\cdots((F_1\ltimes F_2) \ltimes F_3) \ltimes \cdots F_{n-2})
\ltimes F_{n-1}$, using the Magnus ordering in each $F_i$. This
semidirect product decomposition is called the {\it Artin combing}
of $P_n$.

\begin{figure}[ht]
\centerline{\includegraphics{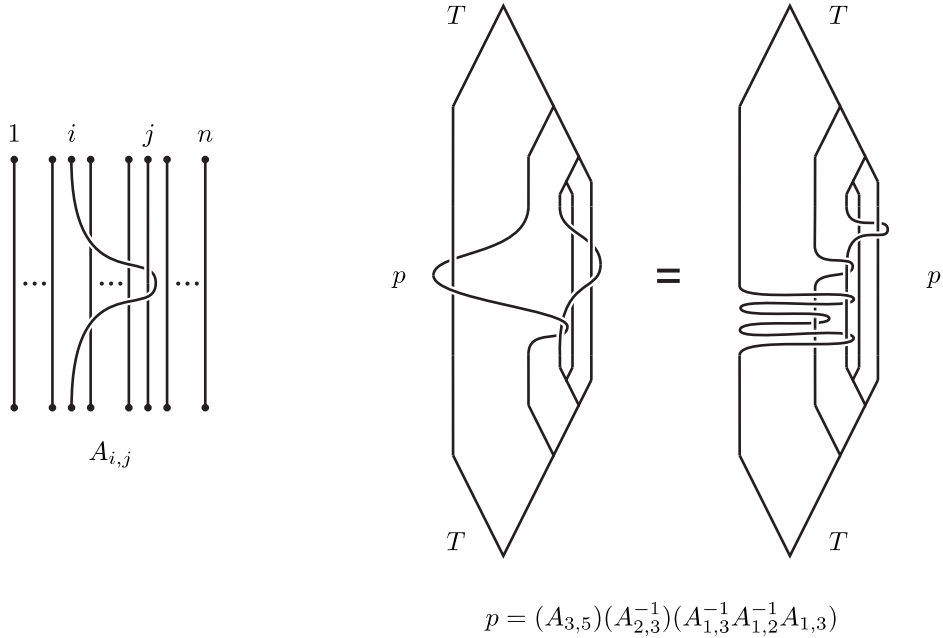}} \caption{The
generators $A_{i,j}$, and the combing of the pure braid $p$.}
\label{puregen_combing}
\end{figure}

It will be convenient to introduce a free generating set for each
$F_i$ in the above semidirect product decomposition of $P_n$. Notice
that each element in $F_i$ corresponds to a loop made by the strand
$n-i$ crossing the strands $n-i+1,n-i+2,\ldots,n$. In this way, if
we consider the pure braid $A_{i,j}$ in the left hand side of
Figure~\ref{puregen_combing}, the group $F_i$ is freely generated by
$\{A_{n-i,j};\ j>n-i\}$.  In the right hand side of
Figure~\ref{puregen_combing} we can see the Artin combing of the
pure braid $p$ in Figure~\ref{bfelem_pbvelem}, which is $p=
(A_{3,5}) (A_{2,3}^{-1}) (A_{1,3}^{-1} A_{1,2}^{-1}A_{1,3})$. Notice
that the strands 4 and 5 do not cross, so the first nontrivial
factor in the decomposition of $p$ is $A_{3,5}$, which is a positive
element of $F_2$. Hence $p$ is a positive pure braid.

\section{The main theorem}

We have already described all the tools we need to show the main
results of this paper, namely that $PBV$ and $BF$ are bi-orderable.
Recall that $\displaystyle PBV=  \lim_{\substack{\displaystyle \longrightarrow \\
T}} P_{n,T}$, that is, $PBV$ is the direct limit of an infinite
number of copies of $P_n$ (with distinct values of $n$). We already
know that each $P_n$ has a bi-order, hence $P_{n,T}$ is bi-orderable
for every $T$, where $(T,p,T)\in P_{n,T}$ is said to be positive if
and only if $p$ is a positive pure braid. But it is not clear that
this order is compatible with the above direct system. This is shown
in the following result, which we will use several times later.

\begin{lem}\label{L:order-alpha compatible}
The element $(T,p,T)\in P_{n,T}$ is positive in $P_{n,T}$ if
and only if $\alpha_{n,T,i}(T,p,T)$ is positive in $P_{n+1,T'}$, for any $i=1,\ldots,n$.
\end{lem}

\begin{proof}
Since the order in $P_n$ is determined by the Artin combing, we need to see
how the Artin combing of a pure braid $p\in P_n$ is
transformed when we apply $\alpha_{n,T,i}$ to the element
$(T,p,T)$. Suppose that $p=f_1f_2\cdots f_{n-1}$ is the Artin
combing of $p$, where each $f_i\in F_i=\langle A_{n-i,l};\ l>n-i
\rangle$, and suppose that $f_j$ is the first nontrivial factor in
the above decomposition, for some $j\geq 1$. This means that the last $j$ strands of $p$
form a trivial braid, while the last $j+1$ strands do not. Now
suppose we apply $\alpha_{n,T,i}$ to $(T,p,T)$. The pure braid $p$
is then replaced by $p'$, which is obtained from $p$ by doubling its
$i$-th strand. Let $p'=f'_1f'_2\cdots f'_n$ be the Artin combing of
$p'$ (notice that $p'$ has one more strand, hence one more factor in
its Artin combing). We will distinguish two cases.

If $i\leq n-j$, the last $j+1$ strands of $p$ and $p'$ form the same
braid, only that the indices in $p'$ are shifted by one. Hence
$f'_1,\ldots,f'_{j-1}$ are all trivial, and $f_j'$ is obtained from
$f_j$ by replacing each $A_{n-j,k}$ by $A_{n-j+1,k+1}$. Since $f_j$
and $f_j'$ determine the same element of $F_j$, their Magnus
expansions coincide, hence $p$ is positive if and only if $p'$ is
positive. The geometric meaning of this computation is that, in this
case, after recombing the braid obtained by splitting the $i$-th strand,
the first factor of the combing is the same as before, which is clear
geometrically.

If $i>n-j$, we have doubled one of the last $j$ strands of $p$,
hence $p'$ is a braid whose last $j+1$ strands form a trivial braid,
and then $f'_1,\ldots,f'_j$ are all trivial. Moreover, $f'_{j+1}$ is
obtained from $f_j$ by replacing each $A_{n-j,k}$ either by
$A_{n-j,k}$ (if $k<i$), by $A_{n-j,k+1}$ (if $k>i$), or by
$A_{n-j,k}A_{n-j,k+1}$ (if $k=i$). To simplify notation, if we
denote $A_{n-j,k}$ by $x_k$, then $f_{j+1}'$ is obtained from $f_j$
by the injective homomorphism $$\theta_i: F_n \rightarrow F_{n+1}$$ given by
$$
\theta_i(x_k)=\left\{\begin{array}{ll}
  x_k & \mbox{ if } k<i, \\
  x_kx_{k+1} & \mbox{ if } k=i, \\
  x_{k+1} & \mbox{ if } k>i.
\end{array}
\right.
$$
.

The problem is now reduced to its algebraic setting. The proof will be
finished if we show that $f\in F_n$ is positive if and only if
$\theta_i(f)\in F_{n+1}$ is positive. We will do this by looking at the
deviations of $f$ and $\theta_i(f)$. Let $F_n=G_1\supset G_2\supset
G_3\supset \cdots$ and $F_{n+1}=G'_1\supset G'_2\supset G'_3\supset
\cdots$ be the lower central series of $F_n$ and $F_{n+1}$,
respectively. Suppose that $f\in G_d$ and $f\notin G_{d+1}$, so its
deviation $\delta(f)$ is a form of degree $d$ in $\mathbb Z[[\mathbf
X]]$. It is well-known~\cite{MKS} that $G_d/G_{d+1}$ is generated by
the elements $[[\cdots[[x_{i_1},x_{i_2}],x_{i_3}],\cdots],
x_{i_d}]$, where $i_k\in\{1,\ldots, n\}$ for $k=1,\ldots,d$, and
$[a,b]=a^{-1}b^{-1}ab$. We can assume $i_1\neq i_2$, since otherwise
the bracket is trivial. These generators are called {\it simple
commutators}, and we remark that they {\it do not} form a basis of
the abelian group $G_d/G_{d+1}$, although they do generate the
group. As in~\cite{MKS}, we will denote each of the above simple
commutators by $[x_{i_1},x_{i_2}\ldots, x_{i_d}]$.

It is shown in~\cite{MKS} that given $x,y\in F_n$ then, if
$\delta(x)\delta(y)-\delta(y)\delta(x)\neq 0$, one has
$\delta([x,y])= \delta(x)\delta(y)-\delta(y)\delta(x)$. Hence, since
$\delta(x_k)=X_k$ for $k=1,\ldots,n$, it follows that
$\delta([x_{i_1},x_{i_2}])=X_{i_1}X_{i_2}-X_{i_2}X_{i_1}$, provided
$i_1\neq i_2$. We will show by induction on $d$ that
$$
\delta([x_{i_1},\ldots, x_{i_d}])=
\delta([x_{i_1},\ldots,x_{i_{d-1}}])X_{i_d}-
X_{i_d}\delta([x_{i_1},\ldots,x_{i_{d-1}}]),
$$
that no monomial in the above expression has the form $X_j^{d}$, and
that the coefficient of the smallest monomial is either 1 or $-1$.
We know the claim is true for $d=2$, so suppose it is true for $d-1$
and let $M$ be the smallest monomial of
$\delta([x_{i_1},\ldots,x_{i_{d-1}}])$. Notice that the smallest
monomial of $\delta([x_{i_1},\ldots,x_{i_{d-1}}])X_{i_d}$ is
$MX_{i_d}$, and the smallest monomial of
$X_{i_d}\delta([x_{i_1},\ldots,x_{i_{d-1}}])$ is $X_{i_d}M$. These
two monomials cannot coincide, since this would imply that
$M=X_{i_d}^{d-1}$, which is not true by induction hypothesis. Hence,
either $MX_{i_d}$ or $X_{i_d}M$ is the smallest monomial of
$\delta([x_{i_1},\ldots,x_{i_{d-1}}])X_{i_d}-
X_{i_d}\delta([x_{i_1},\ldots,x_{i_{d-1}}])$. Moreover, its
coefficient will be $\pm 1$, since the same happens for $M$ in
$\delta([x_{i_1},\ldots,x_{i_{d-1}}])$, by induction hypothesis.
This in particular implies that
$$
\delta([x_{i_1},\ldots,x_{i_{d-1}}])X_{i_d}-
X_{i_d}\delta([x_{i_1},\ldots,x_{i_{d-1}}])\neq 0,
$$
and since $\delta(x_{i_d})=X_{i_d}$ it follows that
$\delta([x_{i_1},\ldots, x_{i_d}])$ is equal to the above
expression, as we wanted to show. It is also clear that no monomial
of $\delta([x_{i_1},\ldots, x_{i_d}])$ can be a power of a variable,
since this would imply that the same happens for some monomial of
$\delta([x_{i_1},\ldots,x_{i_{d-1}}])$. We have then shown the
claim.

A particular consequence of the above claim is that
\begin{equation}\label{E:delta}
\delta([x_{i_1},x_{i_2}\ldots, x_{i_d}])= \sum_{\sigma\in
\Sigma}{\varepsilon_\sigma X_{i_{\sigma(1)}}X_{i_{\sigma(2)}}\cdots
X_{i_{\sigma(d)}}},
\end{equation}
where $\Sigma$ is a certain subset of the symmetric group $S_d$
(which only depends on $d$) and $\varepsilon_\sigma=\pm 1$. Hence,
every monomial in $\delta([x_{i_1},x_{i_2}\ldots, x_{i_d}])$
consists of a permutation of the variables. Notice that if there are
some repeated variables in $X_{i_1},\ldots, X_{i_d}$, then some
monomials in the above expression may coincide, so the coefficients
in the form $\delta([x_{i_1},x_{i_2}\ldots, x_{i_d}])$ may be
distinct from $\pm 1$.

Since we will compare simple commutators in $F_n$ and in $F_{n+1}$,
we will need the following concepts. Given a simple commutator
$[x_{i_1},\ldots,x_{i_d}]\in F_n$ and a simple commutator
$[x_{j_1},\ldots,x_{j_d}]\in F_{n+1}$, we will say that the latter
is an {\em $i$-successor} of the former if it is obtained from it by
replacing each $x_k$ by $x_{k+1}$, if $k>i$,  and each $x_i$ by either
$x_i$ or $x_{i+1}$. Notice that if $m$ is the number of appearances
of $x_i$ in $[x_{i_1},\ldots,x_{i_d}]$, then this commutator has
$2^m$ $i$-successors.  In the same way, if $M$ is a monomial in the
variables $X_1,\ldots,X_n$, and $M'$ is a monomial in the variables
$X_1,\ldots, X_{n+1}$, we will say that $M'$ is an $i$-successor of $M$
if it is obtained from it by replacing each $X_k$ ($k>i$) by
$X_{k+1}$ and each $X_i$ by either $X_i$ or $X_{i+1}$. As above, if
$m$ is the number of appearances of $X_i$ in $M$, then $M$ has $2^m$
$i$-successors. Notice also that among those $2^m$ $i$-successors, there is
only one which does not involve the variable $X_{i+1}$, which we
will call the {\it minimal $i$-successor} of $M$, since it is the
smallest one with respect to the total order on monomials defined
above.

Let us define $\xi_i: \mathbb Z[[X_1,\ldots,X_n]] \rightarrow \mathbb
Z[[X_1,\ldots, X_{n+1}]]$ to be the ring homomorphism given by
$$
\xi_i(X_k)=\left\{\begin{array}{ll}
  X_k & \mbox{ if } k<i, \\
  X_k+X_{k+1} & \mbox{ if } k=i, \\
  X_{k+1} & \mbox{ if } k>i.
\end{array}
\right.
$$
Observe that if $M$ is a monomial in the variables $X_1,\ldots, X_n$,
then $\xi_i(M)$ is equal to the sum of all $i$-successors of $M$.

Now suppose we apply $\theta_i$ to $[x_{i_1},\ldots, x_{i_d}]$. This
would replace each $x_{i_k}$ by either $x_{i_k}$ or
$x_{i_k}x_{i_k+1}$ or $x_{i_k+1}$, depending wether $i_k$ is smaller
than, equal to or greater than $i$, respectively. It is shown
in~\cite{MKS} that for every $a_1,\ldots,
a_{k-1},a_k,a_k',a_{k+1},\ldots, a_d\in F_{n+1}$ (actually in any
group), one has
$[a_1,\ldots,a_ka_k',\ldots,a_d]=[a_1,\ldots,a_k,\ldots,a_d]\:
[a_1,\ldots,a_k',\ldots,a_d]$ (mod $G'_{d+1}$). Hence, if we define
$\mathcal S$ to be the set of all $i$-successors of
$[x_{i_1},\ldots,x_{i_d}]$, one has
$$
   \theta_i([x_{i_1},\ldots,x_{i_d}])=\prod_{[x_{j_1},\ldots,x_{j_d}]\in\mathcal S}{[x_{j_1},\ldots,x_{j_d}]} \quad (\mbox{mod } G'_{d+1}).
$$

Now, if we apply $\delta$ to $\theta_i([x_{i_1},\ldots,x_{i_d}])$,
Theorem~\ref{T:deviation} tells us that
\begin{equation}\label{E:delta theta}
\delta(\theta_i([x_{i_1},\ldots,x_{i_d}]))=
\sum_{[x_{j_1},\ldots,x_{j_d}]\in \mathcal
S}{\delta([x_{j_1},\ldots,x_{j_d}])}.
\end{equation}
Let us see that the above form can be rewritten as follows:
$$
\delta(\theta_i([x_{i_1},\ldots,x_{i_d}]))=
\xi_i(\delta([x_{i_1},\ldots,x_{i_d}])).
$$
Indeed, by~(\ref{E:delta}) and~(\ref{E:delta theta}) one has
$$
\delta(\theta_i([x_{i_1},\ldots,x_{i_d}]))=
\sum_{[x_{j_1},\ldots,x_{j_d}]\in \mathcal S}{\left(\sum_{\sigma\in
\Sigma}{\varepsilon_\sigma X_{j_{\sigma(1)}}X_{j_{\sigma(2)}}}\cdots
X_{j_{\sigma(d)}}\right)}
$$
$$
= \sum_{\sigma\in \Sigma}{\varepsilon_\sigma\left(
\sum_{[x_{j_1},\ldots,x_{j_d}]\in \mathcal S}{
X_{j_{\sigma(1)}}X_{j_{\sigma(2)}}\cdots X_{j_{\sigma(d)}}}\right)}
$$
$$
 = \sum_{\sigma\in \Sigma}{\varepsilon_\sigma\:
\xi_i\left( X_{i_{\sigma(1)}}X_{i_{\sigma(2)}}\cdots
X_{i_{\sigma(d)}}\right)}
$$
$$
= \xi_i\left(\sum_{\sigma\in \Sigma}{\varepsilon_\sigma\:
X_{i_{\sigma(1)}}X_{i_{\sigma(2)}}\cdots X_{i_{\sigma(d)}}}\right)
$$
$$
= \xi_i(\delta([x_{i_1},\ldots,x_{i_d}])).
$$

Finally, suppose that $f\in F_n$ belongs to $G_d$ but not to
$G_{d+1}$. Since the set of simple commutators generates
$G_d/G_{d+1}$, one has $f=c_1^{e_1}\cdots c_t^{e_t}$ (mod
$G_{d+1}$), where each $c_j$ is a simple commutator of order $d$. By
Theorem~\ref{T:deviation} and by the fact that $\theta_i$ and $\xi_i$
are homomorphisms, it follows that $\delta(\theta_i(f))=\xi_i(\delta(f))$,
since the same equality is true for the generators $c_j$.


Therefore, every monomial in $\delta(\theta_i(f))$ with nontrivial
coefficient is an $i$-successor of a unique monomial of $\delta(f)$, and
furthermore their coefficients coincide. This means that the
smallest monomial in $\delta(\theta_i(f))$ is precisely the minimal
$i$-successor of the smallest monomial of $\delta(f)$. Since their
coefficients coincide, it follows that $f$ is positive if and only
if so is $\theta_i(f)$, as we wanted to show. The proof of Lemma~\ref{L:order-alpha compatible} is finished.\end{proof}

From Lemma~\ref{L:order-alpha compatible} we see that doubling the
$i$-th strand preserves the order of the pure braid group, so
$\alpha_{n,T,i}$ preserves the order of $P_{n,T}$. It follows that
if $(T,p,T)$ and $(T',p',T')$ are two representatives of the same
element in $PBV$, then $p$ is a positive pure braid if and only if
so is $p'$. This allows to define an ordering in $PBV$, just by
saying that $(T,p,T)$ is positive if and only if $p$ is a positive
pure braid. By Lemma~\ref{L:order-alpha compatible} this is well
defined, and by the following result, it is a bi-order.

\begin{cor}\label{C:PBV bi-orderable}
The group $PBV$ is bi-orderable.
\end{cor}

\begin{proof}
Since $P_n$ is bi-orderable for every $n$, and $P_n$ is isomorphic
to $P_{n,T}$, we can define $\mathcal P_{n,T}$ to be the semigroup
of positive elements of $P_{n,T}$, that is $\mathcal
P_{n,T}=\{(T,p,T)\in P_{n,T};\ 1<p \}$. Now define $\displaystyle
\mathcal P= \lim_{\substack{\displaystyle \longrightarrow \\ T}}
{\mathcal P_{n,T}} \subset PBV$, that is,
$\mathcal P$ is the set of elements in $PBV$ having one
representative $(T,p,T)$ such that $p$ is positive. By
Lemma~\ref{L:order-alpha compatible}, $\mathcal P$ is also the set
of elements in $PBV$ {\it all} of whose representatives have the
form $(T,p,T)$ with $p$ a positive braid. We will show that this set
defines a bi-order in $PBV$. We must show that $\mathcal P$ is a
semigroup, that $PBV=\mathcal P \sqcup \{1\}\sqcup \mathcal P$, and
that $\mathcal P$ is invariant under conjugation.

We see that $\mathcal P$ is a semigroup, since two elements
$(T_1,p_1,T_1),(T_2,p_2,T_2)\in \mathcal P$ can be multiplied by
finding suitable representatives with matching trees, say
$(T,p_1',T),(T,p_2',T)$, and their product will be $(T,
p_1'p_2',T)$. Since $p_1$ and $p_2$ are positive pure braids, from
Lemma~\ref{L:order-alpha compatible}, the same happens to $p_1'$ and
$p_2'$, hence $p_1'p_2'$ is also positive, and so is $(T,
p_1'p_2',T)\in \mathcal P$. On the other hand, since the inverse of
$(T,p,T)\in PBV$ is $(T,p^{-1},T)$, it follows immediately that $PBV
=\mathcal P \sqcup \{1\} \sqcup \mathcal P^{-1}$. Finally, the
action of conjugating an element $(T_1,p_1,T_1)\in \mathcal P$ by
another element $(T_2,p_2,T_2)\in PBV$ is done by obtaining suitable
representatives $(T,p_1',T)\in \mathcal P$ and $(T,p_2',T)\in PBV$,
and the result is $(T,p_2'p_1'(p_2')^{-1},T)$. Since the set of
positive pure braids is invariant under conjugation, it follows that
the resulting pure braid is positive, so $\mathcal P$ is invariant
under conjugation. Therefore $\mathcal P$ defines a bi-order on
$PBV$, as we wanted to show.
\end{proof}

\medskip

\begin{thm}\label{T:BF bi-orderable}
The group $BF$ is bi-orderable.
\end{thm}

\begin{proof}
This is a consequence of the above corollary, together with
Proposition~\ref{P:orders and extensions}, since by definition of
$BF$ one has a short exact sequence
$$
 1 \longrightarrow PBV \longrightarrow BF
 \stackrel{\widetilde\rho}{\longrightarrow} F \longrightarrow 1.
$$
Moreover, this sequence splits, a section of $\widetilde\rho$ being
the map that sends $(T_-,T_+)\in F$ to $(T_-,1,T_+)\in BF$. Hence
$BF=F \ltimes PBV$.  By Proposition~\ref{P:orders and extensions},
we just need to show that the action of $F$ on $PBV$ preserves the
order of $PBV$ defined in Corollary~\ref{C:PBV bi-orderable}. The
action determined by an element $(T_-,T_+)\in F$ on an element
$(T,p,T)\in PBV$ is given by conjugating the latter by
$(T_-,1,T_+)$, that is, if we choose representatives
$(T_-',1,T')=(T_-,1,T_+)$ and $(T',p',T')=(T,p,T)$, then the conjugate is $(T_-',1,T')(T',p',T')(T',1,T_-')=(T_-',p',T_-')$.
Hence the pure braid $p$ is replaced by $p'$. Since
$(T,p,T)=(T',p',T')$, by Lemma~\ref{L:order-alpha compatible},
it follows that $p$ is positive if and only if
so is $p'$. Therefore, the action of $F$ on $PBV$ preserves the
order, so $BF$ is bi-orderable as we wanted to show, an explicit
bi-order being the lexicographic order in $F\ltimes PBV$, where $F$
is given the usual bi-order and $PBV$ is
bi-ordered as in Corollary~\ref{C:PBV bi-orderable}.
\end{proof}

Notice that the order in $BF$ can be easily described. An element
$(T_-,p,T_+)$ is positive if $(T_-,T_+)\in F$ is positive (its first
slope different from 1 is greater than one), or if $(T_-,T_+)$ is
trivial ($T_-=T_+$) and $p$ is a positive pure braid.


\begin{thebibliography}{99}

\bibitem{Artin1925}  E. Artin,  {\it Theorie der Z\"opfe},  Abh. Math. Sem. Hamburg, {\bf 4} (1925), 47-72.

\bibitem{Artin} E. Artin, {\it Theory of braids},  Annals of Math. {\bf 48} (1946) 101-126.

\bibitem{Birman} J. S. Birman, {\it Braids, Links, and Mapping Class Groups.} Annals of
Math.~Studies, 82, Princeton University Press 1975.

\bibitem{Birman-Brendle} J. S. Birman, T.E. Brendle, {\it Braids: A Survey}. www.arxiv.org/math.GR/0409205

\bibitem{BBCS} T. Brady, J. Burillo, S. Cleary and M.
Stein, {\it Pure braid subgroups of braided Thompson's groups,}
www.arxiv.org/math.GR/0603548

\bibitem{brin} M. G. Brin, {\it The Algebra of Strand Splitting. I. A Braided Version of Thompson's Group
V}. www.arxiv.org/math.GR/0406042

\bibitem{BG} R. Geoghegan, K. S. Brown, {\it An infinite-dimensional torsion-free $FP_\infty$ group}, Invent. Math. {\bf 77} (1984) pp- 367--381.

\bibitem{CFP} J. W. Cannon, W. J. Floyd and W. R. Parry, {\it Introductory notes on Richard Thompson's
groups.} Enseign. Math. (2), {\bf 42} (3-4):215-256, 1996.

\bibitem{dehornoy} P. Dehornoy, {\it From large cardinals to braids via
distributive algebra}, J. Knot Theory Ramifications {\bf 4} (1995),
no. 1, 33-79.

\bibitem{dehornoyBV} P. Dehornoy, {\it The group of parenthesized braids,} Adv. in Math., to
appear.

\bibitem{5authors} R. Fenn, M. T. Greene, D. Rolfsen, C. Rourke and
B. Wiest, {\it Ordering the braid groups}, Pacific J. Math. {\bf 191}
(1999), no. 1, 49-74.

\bibitem{Hansen} V. L. Hansen, {\it Braids and Coverings}, Cambridge University Press, 1989, London Mathematical Society Student Texts, 18.

\bibitem{KR}
D. M. Kim and D. Rolfsen, {\it An ordering for groups of pure braids
and fiber-type hyperplane arrangements.} Canad. J. Math. {\bf 55},
no. 4 (2003), 822-838.

\bibitem{MKS} W. Magnus, Karras and Solitar, {\it Combinatorial Group
Theory. Presentations of groups in terms of generators and
relations}, Revised ed., Dover Publications, Inc., New York, 1976.

\bibitem{MR} R. Mura and A. Rhemtulla, {\it Orderable Groups}. Lecture Notes in Pure and Appl. Math. 27,
Marcel Dekker, New York, 1977.



\end{thebibliography}

\end{document}